\documentclass[11pt]{article}
   \usepackage{amsthm}
   \usepackage{amsmath}
   \usepackage{amssymb}
   \usepackage{amsfonts}
   \usepackage{amscd}
   \usepackage[mathscr]{eucal}
   \usepackage{verbatim}
   \usepackage{graphics}
\usepackage{hyperref}  
  
\newcommand{\aaa}{\mathcal{A}}
\newcommand{\bbb}{\mathcal{B}}
\newcommand{\ccc}{\mathcal{C}}

\newcommand{\ggG}{\mathcal{G}}

\newcommand{\baaa}{\bar{\aaa}}
\newcommand{\kkk}{\mathcal{K}}

\newcommand{\bs}{{\bbb}^*}

\newcommand{\ad}{{\aaa}^{\sharp}}

\newcommand{\kd}{{\kkk}^{\sharp}}

\newcommand{\zp}{\mathbb{Z}/p\mathbb{Z}}

\newcommand{\rr}{\mathbb{R}}

\newcommand{\zz}{\mathbb{Z}}
\newcommand{\sip}{\sigma_p}
\newcommand{\Card}{\operatorname{Card}}

\newcommand{\card}{\operatorname{card}}
\newcommand{\intsqdp}{\lfloor\sqrt{2p}\rfloor}

\newtheorem{thm}{Theorem}
\newtheorem*{mlem*}{Main Lemma}
\newtheorem*{thm*}{Theorem}
\newtheorem{prop}[thm]{Proposition}
\newtheorem{notat}[thm]{Notation}
\newcommand{\N}{{\mathbb N}}
\newcommand{\Z}{{\mathbb Z}}
\newtheorem{lem}[thm]{Lemma}
\newtheorem*{cor*}{Corollary}
\newtheorem*{conj*}{Conjecture} 

\newtheorem{conj}[thm]{Conjecture}
\theoremstyle{definition}

\newtheorem{rem}[thm]{Remark}

\newtheorem{exams}[thm]{Examples}

\begin{document}

\title{Large zero-free subsets of $\zp$} \author{Jean-Marc
Deshouillers, Gyan Prakash\\  
} \date{} \maketitle
\begin{abstract}
A finite subset $\aaa$ of an abelian group $\ggG$ is said to be zero-free if
the identity element of $\ggG$  cannot be written as a sum of distinct
elements from $\aaa.$ In this article we study the structure of
zero-free subsets of $\zp$ the cardinality of which is close to largest possible. In
particular, we determine the cardinality of the largest zero-free subset of
$\Z/p\Z$, when $p$ is a sufficiently large prime.
\end{abstract}

\bigskip

For a finite abelian group $(\ggG,~+)$ and a subset $\aaa$ of $\ggG$, we set
$\ad = \{\sum_{b \in \bbb}b~: \bbb \subset~\aaa,~\bbb~\neq~\emptyset~\}$. We say $\aaa$ is {\it zero-free} if $0 \notin \ad$; in
other words $\aaa$ is zero-free if $0$ can not be expressed as a sum of
distinct elements of $\aaa.$\\

In 1964, Erd\H os and Heilbronn \cite{EH} made
the following conjecture, supported by examples showing that the upper bound
they conjectured is, if correct, very close to being best possible.

\begin{conj} \label{Erd-Heil} 
Let $\aaa$ be a subset of $\zp$.
 If $\aaa$ is zero-free, we have $\Card (\aaa) \le  \sqrt{2p}$.
\end{conj}
\smallskip

Up to
recently, the best result concerning zero-free subsets of $\zp$ was that of
Hamidoune and Z\'emor~\cite{ham} who proved in 1996 that their cardinality is
at most $\sqrt{2p} + 5 \ln p$, thus
showing that the constant $\sqrt{2}$ in the above conjecture is sharp.\\

The study of this question has been revived more recently. Freiman and the
first named author introduced a method based on trigonometrical sums which led
to the description of large incomplete subsets~\cite{DF} as well as that of
large zero-free subsets~\cite{D0} of $\zp$. Recall that a subset
$\aaa$ of $G$ is said to be {\it incomplete} if $\ad \cup \{0\}$ is not
equal to $G.$ Szemer\'edi and Van Vu
~\cite{szvv}, as a consequence of their result on long arithmetic progressions
in sumsets, gave structure results for zero-free subsets
leading to the optimal bound for the total number of such subsets of $\zp$. As
it was noticed independently by Nguyen, Szemer\'edi and Van Vu~\cite{nszvv} on
one side and us on the other one, both methods readily lead to a proof of the
Erd\H os-Heilbronn conjecture for zero-free subsets\footnote{Van H. Vu and the
first named author exchanged this information during a private conversation
held in Spring 2006.}.\\

The aim of the present paper is to study the description of rather large
zero-free subsets of $\zp.$ We start by reviewing the present knowledge on zero-free subsets
of $\zp$.
\begin{notat}
 We denote by $\sip$ the canonical homomorphism from $\mathbb{Z}$
onto $\zp$; for an element $a$ in $\zp$, we denote by $\bar{a}$ be the integer
in $(-\frac{p}{2}, \frac{p}{2}]$ such that $a=\sip(\bar{a})$ and let $|a|_p=
|\bar{a}|$. Given a set $\aaa \subset \zp$, we denote by $\bar{\aaa}$ the
set $\{ \bar{a} :  a \in \aaa\}$. For $d \in \zp$, we write $d \cdot \aaa := \{da : a\in \aaa\}$.
Given any real numbers $x, y$ with $x \leq y$, we write $[x, y]_p$ to denote
the set $\sip([x, y] \cap \Z)$. Given a set $\mathcal{B} \subset \Z$ and non negative real numbers $x,
y$, we write $\mathcal{B}(x, y)$ to denote the set $\{ b \in \mathcal{B}: x \leq |b| \leq y \}$
and simply write $\mathcal{B}(x)$ to denote the set $\mathcal{B}(0,x)$.
\end{notat}

\noindent
It is evident that $\mathcal{A} \subset \Z/p\Z$ is zero-free if and only if the set
$(\bar{A})^{\sharp}$ does not contain any multiple of $p$.
This leads to the following examples of zero-free subsets of $\Z/p\Z$.
\begin{exams}\label{exam}
(\textit{i}) Any subset $\mathcal{A}$ of $\Z/p\Z$ which satisfy the properties that $\bar{A}$ is a
subset of
$[1, \frac{p}{2}]$ and 
$\sum_{\bar{a} \in \baaa} |\bar{a}| \leq p -1$
is a zero-free subset of $\Z/p\Z$.\\

\noindent
(\textit{ii}) Given any integer $k$ with $k(k+1)/2 \leq p+1$, the subset $\aaa$ of
$\zp$
with $\aaa = \{-2,1\}_p \cup [3,k]_p$ is a zero-free subset of $\zp$ which has
cardinality equal to $k$.
\end{exams}

Moreover, one readily sees that if a subset $\aaa$ of $\zp$ is zero-free, then it is also the case for the set $s\cdot \aaa$, for any $s$ coprime with $p$.\\

Building on~\cite{DF}, the first named author proved in~\cite{D0} the following result

\begin{thm}\label{dezou}
Let $c>1$, $p$ a sufficiently large prime and $\aaa$ a zero-free subset of $\zp$ with cardinality larger than $c\sqrt{p}$. Then , there exists $d$ coprime with $p$ such that\\
\begin{equation}\label{eqdezou}
\sum_{a\in \aaa} |da|_p < p + O(p^{3/4}\ln p) \, \text{ and } \sum_{a\in \aaa, \overline{da}<0} |da|_p = O(p^{3/4}\ln p),
\end{equation}
where the constants implied in the $O$ symbol depend upon $c,$
\end{thm}
\noindent and built examples showing moreover that none of the above error-terms can be replaced by $o(p^{1/2})$.\\

The error-terms in (\ref{eqdezou}) were reduced to the best possible $O(p^{1/2})$ by Nguyen, Szemer\'edi and Van Vu in~\cite[Theorem 1.9]{nszvv}.\\

The above mentioned paper of Szemer\'edi and Van Vu~\cite{szvv} implicitly contains the following result, formally stated in~\cite{nszvv} as Theorem 2.1.
\begin{thm}\label{nsv21}
Let $\aaa$ be a zero-free subset of $\zp$. Then for some non zero element $d
\in \zp$ the set $d \cdot \aaa$ can be partitioned into two disjoint sets $\aaa '$ and $\aaa ''$, where
\begin{enumerate}
\item $\aaa '$ has negligible cardinality: $|\aaa '|=O(p^{1/2}/\log^2 p)$.
\item We have $\aaa'' \subset [1, p/2]_p \text{ and } \sum_{a'' \in \aaa''} |a''|_p \leq p -1$.
\end{enumerate}
\end{thm}

We first consider the maximal zero-free subsets of $\zp$. The description given in the following theorem is a synthesis of the results established in Sections 1 and 2. 

\begin{thm}\label{largest}
Let $p$ be a sufficiently large  prime and $\aaa$ a zero-free subset of $\zp$ with maximal cardinality. Then 
\begin{equation}\label{cardmax}
\card (\aaa) \text{ is the largest integer } k \text{ such that } k(k+1)/2\le p+1,
\end{equation}
and one may thus write
$\card (\aaa) = \left[\sqrt{2p +9/4} -1/2 \right]= \left[\sqrt{2p}\, \right] - \delta(p), \text{ with } \delta(p) \in \{0, 1\}$. \\
\indent Furthermore, there exists a non-zero element $d$ in $\zp$ such
that the set $d.\aaa$ is the union of two sets $\aaa'$ and $\aaa''$, with
\begin{enumerate}
\item $\aaa' \subset [-2 (1 + \delta(p)), -1]_p \; , \; \; \aaa''
  \subset [1, p/2]_p \; , \; \; \aaa'' \cap (-\aaa')=\emptyset \text{ and }  \card(\aaa') \leq 1+ \delta(p)$,
\item $\sum_{a' \in \aaa'} |a'|_p \leq 2 (1 + \delta(p)) \text{ and } \sum_{a'' \in \aaa''} |a''|_p \leq p -1 + 3\delta(p).$
\end{enumerate}

\end{thm}

The Reader will find a more detailed description of extremal zero-free sets in Section 2. In this Introduction, we limit ourselves to a few remarks and examples.\\

Writing $ \sqrt{2p +9/4} - 1/2 = \sqrt{2p} + \alpha_p -1/2$, we have $\alpha_p=
O(1/\sqrt{p})$. 
One readily sees that $\delta(p)$ takes the values $1$ or $0$ according as the
fractional part of $\sqrt{2p}$ is smaller than $1/2 -\alpha_p$ or larger. Thus the density of the primes $p$ for which the maximal zero-free of $\zp$ subset has cardinality $\left[\sqrt{2p}\right]$ is $1/2$.\\

The sum $\sum_{a'' \in \aaa''} |a''|_p $ can take the values $p+1$ or $p+2$ only in very special cases, namely when one of $p+2, p+3, p+4, p+5, p+6, \text{ or } p+7$ is a value of the polynomial $x(x+1)/2$ at some integral point $x$. The number of such primes $p$ up to $P$ is $O(\sqrt{P})$; the existence of infinitely many such primes is not known and would result from the validity of some standard conjectures, like Schinzel's hypothesis. The set $\aaa = \{-3, 1, 4, 5, 6, \cdots 14, 15\}_{113}$ is an example of a zero-free subset of $\zz/113 \zz$ which satisfies Theorem \ref{largest} with $\card(\aaa)=\intsqdp -1$, $\; \sum_{a'' \in \aaa''} |a''|_p = p+2$ and $p+7=x(x+1)/2$.\\

We now turn our attention to very large zero-free subsets $\aaa$ of $\zp$, i.e. subsets such that $\sqrt{2p}-\card(\aaa) = o(\sqrt{p})$. From now on, we fix a function $\psi$ from $[2, \infty) \text{ to } \rr^+$ which tends to $0$ at $\infty$ and assume that 
\begin{equation}\label{xxl}
e(\aaa) := |\sqrt{2p}-\card(\aaa)| \le \psi(p)\sqrt{p} \; \text{ and } p \text{ is sufficiently large, }
\end{equation}
the term sufficiently large implicitly refereing to the function $\psi$.\\

The following result gives the structure of large zero-free subsets of
$\zp$. It shows that any given large zero-free subset
$\mathcal{A}$ has a dilate,              which is a union of sets
$\mathcal{A}'$ and $\mathcal{A}''$, where $\mathcal{A}''$ is a set closely related to the one given in Example~\ref{exam} $(\textit{i})$ and
the cardinality of $\mathcal{A}'$ is small.

\begin{thm}\label{sl0}
When $p$ is sufficiently large, then given any zero-free subset $\mathcal{A}$ of $\zp$ with $ e(\aaa)$ satisfying~\eqref{xxl},
there exists a non-zero element $d \in \zp$ such that $d\cdot \mathcal{A}$ can be partitioned
into disjoint sets $\aaa'$ and $\aaa''$ with the following properties
\begin{enumerate}
\item The set $\bar{\mathcal{A}''} $ is included in $ [1, \frac{p}{2})$ and we have $\sum_{a'' \in \mathcal{A}''} |a''|_p \leq p - 1$. 
\item The set $\bar{\mathcal{A'}} $ is included in $ [-c\, e(\aaa), c\, e(\aaa)]$ for some absolute constant $c$ and the cardinality of $\mathcal{A}'$ is $ O\left(\sqrt{e(\aaa) +2}\; \ln(e(\aaa) + 2)\right),$
\end{enumerate}
where $e(\aaa)$ is defined in~\eqref{xxl}.
\end{thm} 
To prove Theorems~\ref{largest} and~\ref{sl0}, we prove the following proposition.
\begin{prop}\label{mainprop}
Let $p$ be a prime and $\aaa$ a zero-free subset $\mathcal{A}$ of
$\zp$ with $e(\aaa)$ satisfying~\eqref{xxl}. When $p$ is sufficiently large,
there exists a non-zero element $d \in \zp$ such that
\begin{eqnarray}\label{dense-error}
\sum_{a \in \mathcal{A}}|da|_p &\leq& p +
O\left({e(\aaa)}^{3/2}\ln{(e(\aaa)+ 2)}\right),\\
\sum_{a \in \mathcal{A}, \bar{da}<0}|da|_p &=&  O\left({e(\aaa)}^{3/2}\ln{(e(\aaa)+ 2)}\right)
\end{eqnarray}
\end{prop} 
\begin{rem}\label{prop=conj}
Noticing that for any zero-free subset $\mathcal{A}$ of $\Z/p\Z$, the
corresponding set $\bar{\mathcal{A}} \subset \Z$ can contain at most one element from
the set $\{x, -x\}$ for any integer $x$ we have 
$\sum_{\bar{a} \in \bar{\mathcal{A}}}|\bar{a}| \geq \frac{|\bar{\mathcal{A}}| (|\bar{\mathcal{A}}| +1)}{2}$.
Using this, Conjecture~\ref{Erd-Heil} is an immediate corollary of Proposition~\ref{mainprop}. 
\end{rem}
To prove Proposition~\ref{mainprop} we use Theorem~\ref{dezou} and the
following result from~\cite{DF}.

\begin{thm}(\cite[Theorem 2]{DF})\label{api}
Let $I > L > 100$ and $B > 2C\ln L$ be positive integers such that
\[
C^2 > 500L(\ln L)^2 + 2000I \ln L.
\]
Let $\mathcal{B}$ be a set of $B$ integers included in $[-L, L].$ Then there
exist
$d > 0$ and a subset $\mathcal{C}$ of $\mathcal{B}$ with cardinality $C$ such
that
\begin{enumerate}
\item all the elements of $\mathcal{C}$ are divisible by $d$,
\item $\mathcal{C}^*$ contains an arithmetic progression with $I$ terms and
  common difference $d$,
\item at most $C\ln L$ elements of $\mathcal{B}$ are not divisible by $d$.
\end{enumerate}
\end{thm}

\section{Proof of Proposition~\ref{mainprop}}
Let              $p$ be a sufficiently large prime and $\mathcal{A} \subset \zp$ be as given in Proposition~\ref{mainprop}.
From Theorem~\ref{dezou}, there exists a non-zero element $d \in \Z/p\Z$ such that
(\ref{eqdezou}) holds. Without loss of generality, we may indeed assume that $d=1$ or, equivalently, replace $d\cdot \mathcal{A}$ by $\mathcal{A}$. We then get 
\begin{equation}
\sum_{\bar{a} \in \baaa}|\bar{a}| = \sum_{a \in \mathcal{A}} |a|_p \leq p + O(p^{3/4}\ln p).
\end{equation}
We prove Proposition~\ref{mainprop}, by showing that if $\baaa \subset
[-\frac{p}{2}, \frac{p}{2}]$ is as above then we have
\begin{equation}\label{red-z}
\sum_{ \bar{a} \in \bar{\mathcal{A}}}|\bar{a}| \leq p + O\left({e(\aaa)}^{3/2}\ln{(e(\aaa)+ 2)}\right).
\end{equation}
We shall first show how one can deduce~\eqref{red-z} from the following proposition.

\begin{prop}\label{zprop}
Let $p$ be a sufficiently large prime and $\kkk \subset \Z$ such that
$\kd$ does not contain any multiple of $p$. We recall that $\psi$ is a fixed function from $[2, \infty) \text{ to } \rr^+$ which tends to $0$ at $\infty$. Let us suppose that we have
\begin{equation}\label{assumed-ub}
e(\kkk) := |\sqrt{2p}-\card(\kkk)| \le \psi(p) \sqrt{2p} \; \text{ and } \; \sum_{k \in \kkk} |k| \leq p + s(\kkk), \text{ with } 0 \le s(\kkk) \le p^{0.9}.
\end{equation}
Then, we have in fact 
\begin{equation}\label{eqzprop}
\sum_{k \in \kkk} |k| \leq p + O\left(\kappa^{3/2}\ln \kappa\right),
\end{equation}
where $\kappa = s(\kkk)/\sqrt{p} + e(\kkk) + 2$. Moreover we have
$$\min\{\sum_{k \in \kkk, k > 
  0}|k|\, , \, \sum_{k \in \kkk, k <
  0}|k|\, \} = O( \kappa^{3/2}\ln \kappa).$$
\end{prop}

The fact that $\aaa$ is zero-free and Relations (\ref{eqdezou}) and (\ref{xxl}) permit to apply Proposition~\ref{zprop} with $\kkk=\overline{\aaa}$. When $e(\kkk) \ge p^{1/4} $, then (\ref{eqzprop}) directly implies (\ref{red-z}). But, when $e(\kkk) \le p^{1/4}$, we first obtain from (\ref{eqzprop}) the following weaker inequality

\[
 \sum_{ \bar{a} \in \baaa}|\bar{a}| \leq p + O(p^{3/8}\ln p).
\]
As such, it is weaker than (\ref{red-z}) in this case, we may use $s(\kkk)               =p^{3/8}\ln p $, so that $\kappa = e(\kkk)+O(1)$, and a further application of Proposition~\ref{zprop} leads to Relation  (\ref{red-z}).\\

\noindent
To prove Proposition~\ref{zprop} we need a few lemmas.
\begin{lem}\label{strategy}
Let $m \in \mathbb{Z}, \ell \in \mathbb{N} \text{ and let } \mathcal{B}$ be a subset of $[-\ell \, , \, \ell]\cap \mathbb{Z}$. We have 
$$(\{m, \ldots, m+\ell-1\} + \bs)\cap \mathbb{Z} \, = \, (\,[m-\sum_{b\in \bbb, b<0} |b| \, , \, m + \ell - 1 + \sum_{b\in \bbb, b>0} |b|\,]) \cap \mathbb{Z}.$$
\end{lem}
\begin{proof}
We write $k=|\bbb|$ and $\bbb = \{b_1 < b_2 < \ldots < b_h < 0 \le b_{h+1} < \ldots < b_k\}$,
where $h=0$ if all the elements of $\bbb$ are non negative. For $0 \le u \le k$, we define
$$
\beta_u= \begin{cases}
\sum_{i=1}^{h-u}b_i & \text{ if } 0 \le u \le h-1,\\
0 & \text{ if } u=h,\\
\sum_{j=h+1}^{u}b_j & \text{ if } h+1 \le u \le k.
\end{cases}
$$
Simply notice that $\beta_0 = \min\{s : s \in \bs\}, \beta_k=\max\{s : s\in \bs\}$ and that $\{\beta_0 < \ldots <\beta_k\}$ is a subset of $\bs$ such that the difference between two consecutive elements of which is at most $\ell$.
\end{proof}

\begin{lem}\label{induct}
Let $\bbb \subset \Z, c \in \Z, x \in \mathbb{N}, \ell \ge x+1$ be
such that $\bbb(x)^{\sharp}$ contains $[c \, , \, c+\ell]\cap
\Z$. Then, if there exists an integer $y$ in $[x+1 \, , \, \infty)$
such that $\bbb(y)^{\sharp}$ does not contain $([\, c-\sum_{b\in
  \bbb(x+1, y), b<0} |b| \, , \, c + \ell - 1 + \sum_{b\in \bbb(x+1,
  y), b>0} |b|\,]) \cap \mathbb{Z}$, and $z$ is the least such
integer, then  we have
$$z \ge \ell + \sum_{b\in \bbb(x+1, z-1)} |b| + 1.$$
\end{lem}

\begin{proof}
We notice that $\bbb(z)^{\sharp} \supseteqq \bbb(x)^{\sharp} + \bbb(x+1, z)^*$. Lemma~\ref{strategy} implies that if $z$ has the required property, then $z\ge x+2$. Since $z \ge x+2$, the minimal property of $z$ implies that the set $\bbb(z-1)^{\sharp}$ does contain
$$ \mathcal{I}= ([\, c-\sum_{b\in \bbb(x+1, z-1), b<0} |b| \, , \, c + \ell - 1 + \sum_{b\in \bbb(x+1, z-1), b>0} |b|\,]) \cap \mathbb{Z}.$$
By our assumption, the set $\mathcal{I}\cup \bigcup_{b\in \bbb, |b|=z} (\mathcal{I}+{b})$ is not an interval. This implies (special case of Lemma~\ref{strategy}) that $z \ge \ell + \sum_{b\in \bbb(x+1, z-1)} |b| + 1.$
\end{proof}

\begin{lem}\label{a-a}
Let $\kkk$ be as given in Proposition~\ref{zprop}. Then for any $k \in \kkk$, the
element $-k$ does not belong to $\kkk$.
\end{lem}

\begin{proof}
If claim is not true, then evidently $0 \in \kd$ which is contrary to the 
assumption.
\end{proof}
\begin{lem}\label{dense-int}
We keep the notation of Proposition~\ref{zprop}. For $x \le 0.9 \sqrt{2p}$, the cardinality of
$\kkk(x)$ is   $ x +  O(e(\kkk) + s(\kkk)/\sqrt{p})$.
\end{lem}
\begin{proof}
Lemma~\ref{a-a} immediately implies that
the cardinality of $\kkk(x)$ is at most $x$.
Let us suppose that the cardinality of $\kkk(x)$ is $x - \lambda(x)$. Then using Lemma~\ref{a-a}
we get
\[
\sum_{k \in \kkk} |k| \geq 
\sum_{i = 1}^{x - \lambda(x)}i + \sum_{i = x +1}^{\card(\kkk) + 
\lambda(x)}i
\]
Writing each summand in the second sum on the right hand side of the above 
inequality 
as 
$(i - \lambda(x)) + \lambda(x)$ and then noticing that the number of terms in the
second sum is 
$\card(\kkk) - x$, we get the following inequality
\begin{equation}
\sum_{k \in \kkk} |k| \geq 
\sum_{i = 1}^{\card(\kkk)}i + \lambda(x) 
( \card(\kkk) - x). 
\end{equation}
Since $x \leq 0.9\sqrt{2p}$ and  $\card(\kkk) \geq \sqrt{2p} -
e(\kkk) \geq \sqrt{2p} - 
\psi(p) \sqrt{p}$, the second term
in the right hand side of the above inequality is larger
than $0.05\sqrt{2p}\lambda(x)$, whereas the first term is $p -O(e(\kkk)\sqrt{2p})$. Now comparing the above
inequality with~\eqref{assumed-ub} we obtain
\[
\lambda(x) \leq c(e(\kkk) + s(\kkk)/\sqrt{p}),
\]
for some absolute constant $c$. The lemma readily follows from this fact.
\end{proof}
\begin{lem}\label{a-bdd}
We keep the notation of Proposition~\ref{zprop}. The largest integer $y_0$ belonging to $\kkk \cup -\kkk$ satisfies 
$y_o = O\left(e(\kkk)\sqrt{2p}
+ s(\kkk)\right)$. 
\end{lem}
\begin{proof}
 Using
Lemma~\ref{a-a} we obtain
\begin{equation}\label{eqlargestint}
\sum_{k \in \kkk}|k| \geq \sum_{i = 1}^{\card(\kkk) -1}i + y_0.
\end{equation}
Now the first term on the right hand side of the above inequality is
$p-O(e(\kkk)\sqrt{2p})$. Therefore comparing the above inequality
with~\eqref{assumed-ub}, the assertion follows.
\end{proof}
\begin{lem}\label{C}
We keep the notation of Proposition~\ref{zprop} and let $x$ be a sufficiently
large  integer. Suppose that the cardinality of $\kkk(x)$ is at least
$0.99x$. Then
there exists a subset $\mathcal{C}$ of $A(x)$ with $|\mathcal{C}| =
O(\sqrt{x}\ln x)$ such that $\mathcal{C}^{\sharp}$ contains an arithmetic
progression of length $x$ and common difference $d$, with $d \in \{1,2\}$.
\end{lem}
\begin{proof}
 Applying Theorem~\ref{api} with
$\mathcal{B} = \kkk(x), \ L = x, I=x+1, \ C = \lfloor 100 \sqrt{x}\ln x\rfloor$, we get that
there exists a subset $\mathcal{C}$ of $\kkk(x)$ with $|\mathcal{C}| = O(\sqrt{x}
\ln x)$ such that $\mathcal{C}^{\sharp}$ contains an
arithmetic progression of length $x$ and common difference $d$ dividing at
least $0.8x$ elements of $\kkk(x)$. Since $\kkk(x)$ is contained in an interval 
of length $2x$, we obtain that $d \in \{1, 2\}$.
\end{proof}
\begin{lem}\label{a,a+1}
Let $x$ and $\mathcal{C}$ be as in Lemma~\ref{C}.
Then there exists $k \in \kkk(x) \setminus \mathcal{C}$
such that the element $k + 1$ also belongs to $\kkk(x) \setminus \mathcal{C}$.
\end{lem}
\begin{proof}
 Let $L$ be the set consisting of those elements $l \in [1, x]$ such
that
one of the elements $l$ or $-l$ belongs to the set $\kkk(x) \setminus \mathcal{C}$.
Then $L$ is a set of cardinality at least $0.9x$ contained in an interval of
length $x$. Therefore there exists $l \in L$ such that
$\{l, l+ 1, l+ 2, l + 3, l + 4\} \subset L$. Now by the definition of $L$, for any
$0 \leq i \leq 4$, either $l + i \in \kkk(x) \setminus \mathcal{C}$ or $-(l+i) \in \kkk(x)
\setminus \mathcal{C}$. The lemma follows evidently by showing that there exists $i$ with $0 \leq i \leq 3$ for which
one of the following two
sets, $\{ l + i, l +i + 1\}$ and $\{ -(l+i), -(l+i +1)\}$ is included in $\kkk(x)
\setminus \mathcal{C}$. If not, then replacing $\kkk$ by $-\kkk$ if necessary we have that
$\{-l, l + 1, l+ 3, -l -4\} \subset \kkk(x) \setminus \mathcal{C}$. This would contradict the assumption that
$0$ does not belong to $\kd$. Hence the lemma follows.
\end{proof} 
We are now in a position to prove Proposition~\ref{zprop}.

\begin{proof}[Proof of Proposition~\ref{zprop}]
From Lemma~\ref{dense-int}, there is an integer $x$
which satisfies the assumption of Lemma~\ref{C} and at the same time
$x = O(e(\kkk) + s(\kkk)/\sqrt{p})$. For this choice of $x$,
 let  $\mathcal{C}$ be a subset of
  $\kkk(x)$, as provided by Lemma~\ref{C}. From Lemma~\ref{a,a+1} we obtain a subset $\{k, k+1\}$ of
$\kkk(x)\setminus\mathcal{C}$. Then the set $\mathcal{C}_1 = \mathcal{C} \cup
\{k,k+1\}$ is a subset of $\kkk(x)$ with $\card( \mathcal{C}_1) = O(\sqrt{x}\ln x)$ and
  $(\mathcal{C}_1)^{\sharp}$ contains an interval $[y, y +x)$ of length
    $x$.
 With this
  interval $\mathcal{I}$ applying
  Lemma~\ref{strategy} with $B = \kkk(x) \setminus \mathcal{C}_1$, we obtain
  that $\kkk(x)^{\sharp}$ contains the interval $[y - \sum_{k \in \kkk(x)
    \setminus \mathcal{C}_1, k < 0}|k|, y + x + \sum_{k \in \kkk(x)
    \setminus \mathcal{C}_1, k>0}|k|)$ of length $x + \sum_{k \in \kkk(x)
    \setminus \mathcal{C}_1}|k|$. Then using Lemmas~\ref{induct} and
  \ref{dense-int}, after an elementary calculation, it follows that for some
  positive absolute constant $c_0$, the set $\kkk(p/c_0)$ contains the
  interval $[y - \sum_{k \in \kkk(p/c_0) \setminus \mathcal{C}_1, k < 0}|k|, y + x +  \sum_{k
    \in \kkk(p/c_0)\setminus \mathcal{C}_1, k > 0 }|k|)$ 
  of length $x + \sum_{k \in  \kkk(p/c_0)
    \setminus \mathcal{C}_1}|k|$. Replacing $\kkk$ by $-\kkk$ we
  may assume that $y > 0$. Then since $\kd$ does not contain any
  multiple of $p$ we obtain the following inequalities
 \[
\sum_{k \in \kkk(p/c_0)\setminus \mathcal{C}_1 }|k| \leq p -1
\]
and 
\[
\sum_{k \in \kkk(p/c_0), k< 0}|k| \leq \sum_{c_1 \in
  \mathcal{C}_1 }|c_1| + \sum_{k \in \kkk(p/c_0)\setminus
  \mathcal{C}_1, k < 0}|k| \leq \sum_{c_1 \in \mathcal{C}_1 }|c_1| + y.
\]
From  Lemma~\ref{a-bdd} we have that
  $\kkk(p/c_0) = \kkk$. Moreover it is also
  evident from the construction of $\mathcal{C}_1 $ that $\sum_{c_1
    \in \mathcal{C}_1 }|c_1|
  = O(x^{3/2}\ln x)$. Since $y \in \mathcal{C}_1^{\sharp} $ we have $y \leq
  \sum_{c_1 \in \mathcal{C}_1}|c_1|$. Therefore the assertion follows.
  \end{proof}
\section{Proof of Theorem~\ref{largest}}
Let $p$ be a sufficiently large prime and $\aaa$ a zero-free subset of $\zp$
of the largest cardinality. 
From Proposition~\ref{mainprop} and Remark~\ref{prop=conj}, we have
that
$\card(\aaa) \leq \sqrt{2p}.$ Moreover, since for any prime $p$ the
set
$[1, [\sqrt{2p}]-1]_p$ is an example of a zero-free subset of $\zp$,
it follows that $\card(\aaa) = [\sqrt{2p}] -\delta(p)$ with $\delta(p)
\in \{0,1\}.$ We set $$s(\sqrt{2p}) = \sum_{i=1}^{[\sqrt{2p}]}i =
\frac{[\sqrt{2p}][
\sqrt{2p}+1]}{2}.$$ 
From Example~\ref{exam} (ii), it follows that when $s(\sqrt{2p}) \leq p+1$,
then $\delta(p) = 0$. In this section we shall show that $\delta_p =
0$, only when $s(\sqrt{2p}) \leq p+1.$  

\vspace{2mm}
Using Proposition~\ref{mainprop}, there exists a $d \in (\zp)^*$ such
that replacing $\aaa$ by $d.\aaa$, we have
\begin{equation}\label{sum-prop}
\sum_{a \in \aaa}|a|_p \leq p + O(1).
\end{equation}
Using~\eqref{eqlargestint} with $\kkk = \baaa$ and ~\eqref{sum-prop}, the following lemma is
immediate.
\begin{lem}\label{largestint}
The largest integer $y$ in $\baaa \cup -\baaa$ is $O\left(\sqrt{2p}\right).$ 
\end{lem}

Let $G(\aaa)$ be the collection  of all natural numbers $g$ which satisfy the
property that none of the integers $g$ and $-g$ belong to the set
$\bar{\aaa}$, where $\bar{\aaa}$ is the subset of integers as defined earlier.
For the brevity of notation we shall write $G$ to denote the set $G(\aaa)$.
Let $G = \{g_0 < g_1< g_2<.....\}$.

From Lemma~\ref{dense-int} we obtain
that the cardinality of $G(x)$ is $O(1)$ for any $x \leq 0.9\sqrt{2p}.$  The
arguments identical to those used in the proof of Lemma~\ref{dense-int} in fact
leads to the following lemma.
\begin{lem}\label{except-one}
The set $\baaa \cup (-\baaa)$ contains all the integers in $[1, \sqrt{2p}/5]$ with
at most $\delta(p)$ exception.
\end{lem}
\begin{proof}
The lemma is equivalent to showing that in case $\card(\aaa) =
[\sqrt{2p}]$, then $g_0 > \sqrt{2p}/5$, whereas in case $\card(\aaa) =
[\sqrt{2p}]-1$, then $g_1 > \sqrt{2p}/5.$ Suppose that this is not
true. Then if $\delta(p) = 0$, we have
\[
\sum_{ a \in \aaa}|a|_p = \sum_{\bar{a} \in \bar{\aaa}}|\bar{a}|\geq
\sum_{i = 1}^{g_0 - 1} i + \sum_{i = g_0 +
  1}^{\Card{\aaa} +1}i \geq \frac{[\sqrt{2p}+1][\sqrt{2p}+2]}{2} -\sqrt{2p}/5,
\]
whereas in case $\delta(p) =1$, we have
\[
\sum_{ a \in \aaa}|a|_p  \geq
\sum_{i = 1}^{g_0 - 1} i + \sum_{i = g_0 + 1}^{g_1 - 1}i + \sum_{i = g_1 +
  1}^{\Card{\aaa} +2}i \geq \frac{[\sqrt{2p}+1][\sqrt{2p}+2]}{2}
-2\sqrt{2p}/5.
\]
Using the facts that $[\sqrt{2p}] \geq \sqrt{2p}-1$ and for any
integer $i$, we have $[\sqrt{2p}+i] = [\sqrt{2p}] + i$, it follows that
either of these inequalities are contrary to~\eqref{sum-prop}. Hence
the lemma follows.
\end{proof}

\begin{table}

$$
\begin{array}{|c|c|c|c|c|}\hline
\multicolumn{5}{|c|}{\text{ Structure of } \mathcal{A} \text{ when }
  |\mathcal{A}| = \lfloor \sqrt{2 p}\, \rfloor - \delta_p, \text{ with }
  \delta_p \in \{0,1\}}\\ \hline
\multicolumn{5}{|c|}{\text{ Here } s'' = s''(\aaa) = \sum_{\bar{a} \in \baaa,
    \bar{a}>0}\bar{a}.}\\ \hline
\sharp  & \{a\in \mathcal{A} / |a|_p\leq 4 \}& \delta_p & g_0 & (\baaa)^{\sharp} \\ \hline
1 & \{1, 2, 3, 4\} &  & \ge 5 & [1,s''] \\ \hline
2 & \{-1, 2, 3, 4\} & & \ge 5&\{-1\} \cup [1,s'']  \\ \hline
3 & \{1, -2, 3, 4\} &  &\ge 5 &\{-2,-1\} \cup [1,s''] \\ \hline
4 & \{1, 2, 3\} &               1 &              4&[1,s'']      \\ \hline
5 & \{-1,2,3\} & 1 & 4 &\{-1\}\cup [1,s'']              \\ \hline
6 & \{1,-2,3\} & 1 & 4&\{-1,-2\} \cup [1,s''] \\ \hline
7 & \{-1,2,-3\} & 1 & 4&[-4,-1] \cup [1,s'']\\ \hline
8 & \{1,2, 4\} & 1 & 3&[1, s'']\\ \hline
9 & \{-1,2,4\} & 1 & 3&\{-1\}\cup  [1,s''] \\ \hline
10 & \{1,-2,4\} & 1 & 3&\{-2,-1\} \cup [1,s'']\\ \hline
11 & \{1,2,-4\} & 1 & 3&[-4,-1] \cup [1,s'']\\ \hline
12 &\{-1,-2,4\}& 1& 3&[-3,-1]_p \cup [1,s'']                             \\ \hline
13 & \{1,3,4\} & 1 & 2 & [1,s'']\setminus \{2,s''-2\} \\ \hline
14 & \{-1,3,4\} & 1 & 2 &\{-1\} \cup [1,s'']\setminus \{1,s''-2\}\\ \hline
15 & \{1,-3, 4\} & 1 & 2 &\{-3,-2\}\cup [1,s'']\setminus \{s''-2\}\\ \hline
16 & \{2,3,4\} & 1 & 1 &[2, s'']\setminus \{s''-1\}\\ \hline
17 & \{-2,3,4\} & 1 & 1&\{-2,-1\}\cup[1,s'']\setminus \{s''-1\}
\\ \hline
18 & \{2,-3,4\}& 1 & 1 &\{-3,-1\}\cup[1,s'']\setminus \{s''-1\} \\ \hline 
19 &\{2,3,-4\}& 1 & 1&\{-4,-2,-1\}\cup [1,s'']\setminus \{s''-1\} \\ \hline
\end{array}
$$
\caption{Subset sum of a largest zero-free set}
\end{table}

Now we determine all the possible structure of $\bar{\aaa}(\sqrt{2p}/5)$,
first under the assumption that $g_0 \geq 5$.
\begin{lem} \label{5+}When $g_0 \geq 5$, then  replacing $\aaa$ by $-\aaa$, if
necessary, the set $\baaa$ contains the whole interval $[5,
\sqrt{2p}/5)$ with at most $\delta(p)$ exception and $\baaa(4)$ is
equal to one of the three sets described in  the the first three rows of the
second column of Table 1.
\end{lem}
\begin{proof}
Since we have assumed that $g_0 \geq 5$, replacing $\aaa$ by $-\aaa$, if
necessary, we may assume that $3 \in \baaa$. Then the set $\baaa(3)$ is equal
to one of the following four sets, $\{1,2,3\}, \{-1,2,3\}, \{1,-2,3\},
\{-1,-2,3\}$. Since $\aaa$ is zero-free, among these four possibilities, the
last one cannot occur. We verify that in all the other three possible cases
the following always hold

$$ \{1,2,3,4\} \subset \baaa(3)^{\sharp}.$$

This implies that the set $\baaa(4)$ is equal to one of the three sets described in
the second column of the first three rows of Table 1; that is, the set
$\baaa(4)$ is equal to one of the following three sets $\{1,2,3,4\},
\{-1,2,3,4\}, \{1,-2,3,4\}.$ We claim that there does not exist any integer 
$z \in [5, \sqrt{2p}/5]$ with $-z \in \baaa$. The lemma follows immediately
using this claim and Lemma~\ref{except-one}. To verify the claim, suppose that the
claim is not true and $z_0$ is the least integer which violates the claim.
Then since $\{1,2,3,4,5,6\}$ is always a subset of $\baaa(4)^{\sharp}$, we have
that $z_0$ is at least $6$. Now if $z_0 \ne g_0 + 1$, then we have $z_0 - 1
\in \baaa$ and thus $z_0 \in \baaa(3)^{\sharp} + z_0 -1 \subset
\baaa(z_0-1)^{\sharp}$. Since $\aaa$ is zero-free, this implies that $-z_0$ can
not belong to the set $\baaa$ which contradicts the assumption that $z_0$ is
the least integer violating the claim. Thus if the claim is not true then 
$z_0 =g_0 +1$. But in this case $z_0 -2 \in \baaa$ and thus 
$z_0 \in \baaa(3)^{\sharp} + z_0 -2 \subset \baaa(z_0-2)^{\sharp}$. This implies that
$-z_0$ cannot belong to $\baaa$. Hence the claim and thus the lemma hold.
\end{proof}
\begin{lem}\label{subsetsum5}
Let $\aaa$ be as in Lemma~\ref{5+}, $\aaa' = \aaa \cap [\frac{-p}{2},
-1]_p$
and $\aaa'' = \aaa \cap [1, \frac{p}{2}]_p$. Then we have $\aaa'
\subset [-2, -1]_p.$ Moreover the set $(\aaa)^{\sharp}$ contains the
interval $[1, s'']_p$ with $s'' = \sum_{a'' \in \aaa''} |a''|_p$ and
is equal to one of the sets described in the fifth column of the first three
rows of Table 1, the three possibilities corresponding to three
possible structures for $\baaa(4).$  We have 
$$s'' \leq p-1.$$
\end{lem}
\begin{proof}
For any integer $z$ we set
$$s''(z) =\sum_{\bar{a''} \in
  \baaa''(z)} \bar{a''}.$$
We claim that there is an absolute constant $c$ such that for any
integer $z$
with $5
\leq z
\leq \frac{p}{c}$,  the set 
$\left(\baaa(z)\right)^{\sharp}$ contains the interval $[1, s''(z)]$. The claim is easily verified with $ z=5.$
Suppose the claim is not true
and $z_0$ is the least integer violating the claim. Since
using the previous lemma we always have $s''(5) \geq 5 +1=6$, we 
apply Lemma~\ref{induct} with $x=5$ and obtain  the following
inequality.
\begin{equation}\label{lby0}
z_0 \geq s''(z_0-1) +1.
\end{equation}
Using the previous lemma, for any integer $y$ with $y
\in [6, \sqrt{2p}/5],$ we have 
$$s''(y) = 
\frac{y(y+1)}{2} -\sum_{a'\in \baaa'(4)}|a'|  -\epsilon,$$
where $\epsilon = 0$ if $y\leq g_0$ and  $\epsilon = g_0$ if $y > g_0.$
Using this it follows that ~\eqref{lby0} cannot hold with $z_0 \leq
\sqrt{2p}/5$. Therefore we have
 $$z_0 \geq s''(\sqrt{2p}/5) \geq \frac{p}{c},$$ 
where $c$ is an absolute constant. Hence the claim follows. Using the claim and Lemma~\ref{largestint}, it follows that the set $(\aaa)^{\sharp}$
contains the interval $[1, s'']_p.$ Since $\aaa$ is zero-free, it follows that
 $$s'' \leq p-1.$$ 
 Since
$\left(\baaa(\sqrt{2p}/5)\right)^{\sharp}$ contains the interval $[1,
s''(\sqrt{2p}/5)]$, it follows that there is no integer $y \in
\left[-\frac{p}{c}, -\frac{\sqrt{2p}}{5}\right]$ with $y \in \baaa'$.
Using the
previous lemma and Lemma~\ref{largestint}, it follows that $\baaa' = \baaa \cap [-4,-1] \subset [-2,-1].$ 
 Using this it may be easily
verified that the set $(\aaa)^{\sharp}$ is equal to one of the sets described in
the fifth column of the first three rows of Table 1. Hence the lemma follows.
\end{proof}
\begin{thm}\label{card}
Let $p$ be a sufficiently large prime and $\aaa$ a zero-free subset of
$\zp$ of the largest cardinality. Then $\card(\aaa) =
\left[\sqrt{2p}\right] - \delta(p)$, where $\delta(p) = 0$ if
$s(\sqrt{2p}) \leq p+1$ and is equal to $1$ otherwise. In other words,
$\card(\aaa)$ is the largest integer $k$ with the property that
$\frac{k(k+1)}{2}\leq p+1$; that is $\card(\aaa) = \left[\sqrt{2p
    +9/4} -1/2\right].$
\end{thm}
\begin{proof}
From the remarks made in the beginning of this section, it follows that  $\card(\aaa) = [\sqrt{2p}] - \delta(p)$ with $\delta (p) \in \{0,1\}$.
If $s(\sqrt{2p}) \leq p+1$, then the set $\{-2,-1\}_p \cup \left[3,
  [\sqrt{2p}]\right]_p$ is an example of a zero-free subset of $\zp$
and since
$\aaa$ is a largest zero-free subset, we have $\delta(p) = 0$, in this
case. 
Now in case $\delta(p) = 0$, then from the 
remarks made in the beginning of this section there is a $d \in (\zp)^{*}$, such that replacing $\aaa$ by $d.\aaa$, the inequality~\eqref{sum-prop} holds with $d=1$. Using Lemma~\ref{except-one}, it also follows that $g_0 \geq \sqrt{2p}/5 \geq 5.$
Therefore  it follows that 
replacing $\aaa$ by $-\aaa$,  if necessary, the set $\aaa$ is as in
Lemma~\ref{subsetsum5}.  Since $\delta(p) =0$, we also have that
$$s(\sqrt{2p}) \leq s'' + \sum_{a' \in \aaa'} |a'| \leq s'' + 2,$$
where $s''$ is as in the Lemma~\ref{subsetsum5} and is at most $p-1$. 
 Thus $s(\sqrt{2p}) \leq p+1.$
Hence the theorem follows.  
\end{proof}

\begin{lem}\label{fix4}
Let $\aaa$ be a largest zero-free subset of $\zp$ which
satisfy~\eqref{sum-prop}. When $g_0 \leq 4$, then, replacing $\aaa$ by $-\aaa$ if necessary, the set
$\baaa$ contains the whole interval $[5, \sqrt{2p}/5]$.
\end{lem}
\begin{proof}
Since $g_0 \leq 4$, then using Lemma~\ref{except-one}, for any
integer $z \geq 5$ either $z$ or $-z$ belongs to $\baaa$. Replacing
$\aaa$ by $-\aaa$, if necessary, we may assume that
 the integer $5$ belongs to the set
$\baaa$. If the statement of the lemma is not true then there is an integer $z \in [6,
\sqrt{2p}/5]$ with $-z \in \baaa$. Let $z_0$ be the least among such
integers. Then since $-z_0$ belongs to $\baaa$ and $\aaa$ is zero-free, it
follows that $z_0 -5$ does not belong to the set $\baaa$. From the definition
of $z_0$ it follows that $z_0 -5 \leq 4$ and thus $z_0 \leq 9$. In other
words,
 $z_0 \in \{6,7,8,9\}.$ On the other hand we shall show that  $z_0$ cannot be equal to any of this four possible integers.\\

Case 1: If $z_0 = 9$. In this case we have  $\{5,6,7,8,-9\} \subset \baaa.$
 Since $2+7 -9 = 6+5 -9-2 = 3+6 -9 = 7+5 -3-9 =0$, it follows that none of the
 integers in the set $\{2,3,-2,-3\}$ belongs to the set $\baaa$. This
 is in contradiction to Lemma~\ref{except-one}. Thus $z_0$ cannot be equal to
 $9$.\\

Case 2: If $z_0 = 8$. In this case we have that $\{5,6,7,-8\} \subset \baaa$.
Since we have $3+5-8 = -3 +6+5-8 = 1 +7-8 = -4 +7+5-8 = 0$, none of the
integers from the set $\{1,2,3,-3,-4\}$ belongs to $\baaa$. From
Lemma~\ref{except-one}, it follows that $\{-1,-2,4,5,6,7,-8\} \subset \baaa$.
Since $-1+5+4-8=0$, this is in contradiction to the fact that $\aaa$ is
zero-free. Therefore $z_0$ cannot be equal to $8$.\\

Case 3: If $z_0 = 7$. In this case we have  $\{5,6,-7\} \subset \baaa.$
Since $-4 +5+6-7 = 2+5-7 = 1+6-7 =0$, it follows that none of the integers
from the set $\{1,2,-4\}$ belongs to $\baaa$. Now if $4 \in \baaa$, in other
words if $\{4,5,6,-7\} \subset \baaa$, then since we have $-2 +5+4-7 =-3+4+6-7= 0$, it
follows that there is no integer in $\{2,-2,-3\}$ which belongs to $\baaa$. Therefore we have $g_0
=2$ and using Lemma~\ref{except-one}, the set $\{-1,3,4,5,6,-7\}$ is included in
$\baaa$. Since $-1+3+5-7 = 0$, this is in contradiction to the fact that $\aaa$
is zero-free. Therefore it follows that neither the integer $4$ nor $-4$ can
belong to $\baaa$. Therefore using Lemma~\ref{except-one}, we have
$\{-1,-2,5,6,-7\} \subset \baaa$. Since $3 -1 -2 = -3 -1-2+6 = 0$, this implies
that neither the integer $3$ nor $-3$ can belong to $\baaa$. In other words
none of the integers from the set $\{-3,3,-4,4\}$ can belong to $\baaa$. This
is in contradiction to Lemma~\ref{except-one}. Hence $z_0$ cannot be equal to
$7$.\\

Case 4: If  $z_0 =6$. In this case we have  $\{5,-6\} \subset \baaa$. Since
$1+5-6=0$, it follows that the integer $1$ cannot belong to $\baaa.$ We have two  subcases to discuss in this case, the first one  when $g_0 \ne 1$ and the second one when $g_0 =1.$ .\\

In case $g_0 \ne 1$, then we have $-1 \in \baaa$; that is $\{-1,5,-6\} \subset \baaa$. Since
$-1 -6 + 7 =0$, this implies that $-7 \in \baaa$. This in turn implies that 
$-8 \in \baaa$. Thus we have  $\{-1,5,-6,-7,-8\} \subset \baaa.$ Since
$4 +5 -8-1 = -4 -1 +5 = 0$, it follows that none the integers $4$ nor $-4$
belongs to $\baaa$ and hence $g_0 =4$. Since $3+5-8 = 2+5-7= 0$, it follows
that none of the integers from the set $\{2,3\}$ belongs to $\baaa$.
Hence using Lemma~\ref{except-one} we have 
$\{-2,-3,5\} \subset \baaa.$ Since $\aaa$ is zero-free, this is not
possible. Hence if $z_0 = 6$, then $g_0 =1$.

In case $g_0 =1$, then either $3$ or $-3$
belongs to $\baaa.$

 If $3$ belongs to $\baaa$; that is $\{3,5,-6\} \subset
\baaa$, then since $-2+3+5-6 =0$, it follows that $2 \in \baaa.$ Thus we have
$\{2,3,5,-6\} \subset \baaa$. Since $4+2-6= -4+2+3+5-6=0$, it follows that
none of the integers from the set $\{1,-1,4,-4\}$ belongs to $\baaa$. This
is in contradiction to Lemma~\ref{except-one}. 

In case $-3 \in \baaa$, in other words $\{-3,5,-6\} \subset \baaa.$ Since
$4+5 -3-6 =0$, it follows that $-4 \in \baaa,$ that is $\{-3,-4,5,-6\} \subset
\baaa.$ Since $-2-3+5 = 2-3-4+5 =0$, it follows that none of the integers from
the set $\{1,-1,2,-2\}$ can belong to $\baaa.$. This is in contradiction to
Lemma~\ref{except-one}. 
 
Hence we have shown that $z_0 \notin [6,\sqrt{2p}/5]$ and thus the lemma follows.
\end{proof}
\begin{lem} \label{44} Let $\aaa$ be as in the previous lemma. Then the set $\baaa(4)$ is
equal to one of the sets described in the
second column of the last sixteen rows of Table 1. 
\end{lem}
\begin{proof}
 Let $N$ be
the set of integers $n_i$ which belongs to $[1,4]$ with $-n_i \in \baaa$. Then
it follows using the previous lemma that 
\begin{equation}\label{sumn4}
\sum_{n_i \in N} n_i \leq 4.
\end{equation}
This implies that the cardinality of $N$ is at most $2$.\\

\noindent
When $\card(N) = 2$. It  follows from~\eqref{sumn4} that $N$ is either equal to  $\{1,2\}$ or  is equal to $\{1,3\}$; that is, in this case either
$\{-1,-2\}$ or $\{-1,-3\}$ is a subset of $\baaa$. In case $\{-1,-2\}$ is a
subset of $\baaa$, then since $3-1-2= 0$, it follows that $g_0 = 3$ and
 $\baaa(4)$ is equal
to $\{-1,-2,4\}$. In case $\{-1,-3\}$ is a subset of $\baaa$, then
since $4-1-3=0$, it follows that $g_0 =4$ and $\baaa(4)$ is equal to
$\{-1,-3,2\}.$\\

\noindent
When $\card(N) = 1.$ We have the following four sub-cases to
discuss.
\begin{itemize}

\item When $N= \{1\}$. In this case $\baaa(4)$ can be equal to any of the following
three sets, namely, $\{-1,2,3\}, \{-1,2,4\}, \{-1,3,4\}.$

\item When $N=\{2\}$. In this case $\baaa(4)$ can be equal to any of the following
three sets, namely, $\{-2,1,3\}, \{-2,1,4\}, \{-2,3,4\}.$

\item When $N=\{3\}$. Since $1+2-3=0$, in this case either $g_0$ is equal to $1$ or is
equal to $2$. Moreover the set $\baaa(4)$ is equal to one of the following two
sets, namely, $\{-3,1,4\}, \{-3,2,4\}$.

\item $N= \{4\}$. Since $1+3-4 = 0$, it follows that either $g_0$ is equal to $1$
or is equal to $3.$ In this case $\baaa(4)$ is equal to one of the following
two sets, namely $\{-4,1,2\}, \{-4,2,3\}.$
 \end{itemize}
{\it When $\card(N) = 0$. } In this case $\baaa(4)$ is equal to any one of the
following four sets, namely, $\{1,2,3\}, \{1,2,4\}, \{1,3,4\}, \{2,3,4\}.$
\end{proof}
\begin{lem}\label{subsetsum4}
Let $\aaa$ be as in Lemma~\ref{fix4}, $\aaa' = \aaa \cap
[\frac{-p}{2},-1]_p$
and $\aaa'' = \aaa \cap [1, \frac{p}{2}]_p$. Then we have  
$$ \aaa' \subset [-4, -1]_p, \quad \left[5, \sqrt{2p}-9\right]_p \subset
  \aaa'' \subset \left[1, \sqrt{2p}+8\right]_p.$$
Moreover, the set $\aaa^{\sharp}$ contains $[3, s''-3]_p$ with $s'' =
\sum_{a''\in \aaa''}|a''|_p$ and is equal to
one of the set described in the fifth column of the last sixteen rows of Table
1; the sixteen possibilities correspond to sixteen possibilities for
$\aaa(4)$ as given by Lemma~\ref{44}. We also have $s'' \leq p+2.$ 
\end{lem}
\begin{proof}
For any positive integer $z$ we set
$$s''(z) = \sum_{\bar{a''} \in \baaa''(z)}|\bar{a''}|.$$
We claim that there is an absolute constant $c$ such that for any
integer $z$ with $6 \leq z \leq \frac{p}{c}$, the set
$(\baaa(z))^{\sharp}$ contains the interval $[3, s''(z)-3]$. Suppose the
claim is not true and let $z_0$ be the least integer in $[6, \frac{p}{2}]$
such that $(\baaa(z_0))^{\sharp}$ does not contains the interval $[3, s''(z_0)-3]$.
Since the claim is easily verified when $z=6$, it follows that $z_0
\geq 7.$ Moreover we also verify that the length of the interval $[3,
s''(6)-3]$ is at least $7.$ Therefore
using Lemma~\ref{induct} with $x =6$, it follows that
\begin{equation}
z_0 \geq s''(z_0-1)-4 +1.
\end{equation}
Using
Lemmas~\ref{fix4} and~\ref{44}, it follows that the above
inequality does not hold for any $z_0$ with $z_0 \in [6, \sqrt{2p}/5]$. Therefore we
have
$$z_0 \geq s''(\sqrt{2p}/5-1) -3 \geq \frac{p}{c},$$
where $c$ is an absolute constant. Hence the claim follows. Using Lemma~\ref{largestint}, it follows  that
$(\baaa)^{\sharp}$ contains the interval $[3,s''-3].$  Since $\aaa$ is zero-free, it follows that
 $s'' \leq p+2.$
Since $(\baaa(\sqrt{2p}/5))^{\sharp}$ contains the interval $[3,
\frac{p}{c}]$ it follows there is no $y
\in [-p/c, -\sqrt{2p}/5]$ with $y \in \baaa'.$ Then using Lemma~\ref{fix4}, it follows that $\baaa' = \baaa \cap [-4,-1]\subset
 [-4,-1].$  
Using this, it is easy to verify that the set $(\baaa)^{\sharp}$
 is equal to one of the sets described in the fifth column of the last
 sixteen rows of Table 1. We shall now show that
$$\left[5, \sqrt{2p}-9\right]_p \subset
  \aaa''.$$ 
Since $\baaa' \subset [-4, -1]$, this follows by showing that 
$$g_1 \geq [\sqrt{2p}]-8.$$
For proving this we may assume
that $g_1 \leq \sqrt{2p}$.
Then we observe that the following inequality holds
\begin{equation}\label{1stinclu}
 \sum_{\bar{a}\in \baaa}|a| \geq s(\sqrt{2p}) +\left[\sqrt{2p}\right] +1
-g_0 -g_1.
\end{equation}
The left hand side of the above inequality is equal to $s'' =
\sum_{\bar{a'} \in \baaa'} |\bar{a'}|$ and is thus at most $p+6.$
Moreover using Lemma~\ref{except-one} and Theorem~\ref{card}, we have 
$s(\sqrt{2p}) \geq p+2.$
Using this and rearranging the terms of~\eqref{1stinclu}, we obtain
that $g_1 \geq \sqrt{2p} -8.$
We shall now show that
$$\baaa'' \subset \left[1, [\sqrt{2p}]+8\right].$$
This is equivalent to showing that the largest integer $y \in \baaa$ is
at most $[\sqrt{2p}]+8.$ Now we have the following inequality 
$$\sum_{\bar{a} \in \baaa}|\bar{a}| \geq s(\sqrt{2p}) - g_0-[\sqrt{2p}]+y.$$
Rearranging the  terms of the above inequality we obtain the desired
upper bound for $y.$ Hence the lemma follows.
\end{proof}
\begin{thm}\label{struc-largest}
Let $p$ be  a sufficiently large prime and $\aaa$ be a zero-free
subset of $\zp$ of the largest cardinality. We write $\delta(p)$ to
denote the integer  $\left[\sqrt{2p}\right]-\card(\aaa),$ as in
Theorem~\ref{card}.
Then there exists  $d \in
(\zp)^*$ such that the set $d.\aaa$ is union of sets $\aaa'$ and
$\aaa''$ satisfying the following properties:
\begin{enumerate}
\item $\aaa' \subset [-2 (1 + \delta(p)), -1]_p \; , \; \; \aaa''
  \subset [1, p/2]_p \; , \; \; \aaa'' \cap (-\aaa')=\emptyset \text{ and }  \card(\aaa') \leq 1+ \delta(p)$,
\item the set $\aaa''$ contains the whole interval $[5, \sqrt{2p}/5]_p$
  with at most $\delta(p)$ exception,
\item the set $(-\aaa') \cup \aaa''$ contains the whole interval $[1, 4]_p$, with at most $\delta(p)$ exception,
\item the set $(d.\aaa)^{\sharp}$ contains the whole interval
  $[3,s'']_p$ with at most $\delta(p)$ exception, where $s'' = \sum_{a'' \in \aaa''} |a''|_p,$ 
\item $\sum_{a' \in \aaa'} |a'|_p \leq 2 (1 + \delta(p)) \text{ and } \sum_{a'' \in \aaa''} |a''|_p \leq p -1 + 3\delta(p).$
\end{enumerate}
Further, if $s'' = \sum_{a'' \in \aaa''}
  |a''|_p > p-1$, then we have $s(\sqrt{2p}): =
  \frac{\left[\sqrt{2p}\right]\left[\sqrt{2p}+1\right]}{2} \in [p+2, p+7].$
\end{thm}
\begin{proof}
It is sufficient to show that there exists $d\in (\zp)^*$ such that replacing $\aaa$ by $d.\aaa$, the conclusion of the theorem holds with $d=1.$
From Proposition~\ref{mainprop}, there exists  $d \in (\zp)^*$ such
that replacing $\aaa$ by $d.\aaa$, the inequality~\eqref{sum-prop}
holds. Let $g_0$ be the least positive integer which does not belong
to $\aaa \cup -\aaa.$  When $g_0 \geq 5$, replacing $\aaa$ by $-\aaa$
if necessary,
let $\aaa$ be as in Lemma~\ref{5+}. When $g_0 \leq 4$, then replacing
$\aaa$ by $-\aaa$ if necessary, let $\aaa$ be as in Lemma~\ref{fix4}. For such
$\aaa$, let $\aaa' = \aaa \cap [-\frac{p}{2},-1]_p$ and $\aaa'' = \aaa
\cap [1, \frac{p}{2}]_p.$
Then claims (i)-(v) follow from Lemmas~\ref{except-one},\;
\ref{5+} and \ref{subsetsum5} in case $g_0
\geq 5$ and from Lemmas~\ref{fix4},\;  \ref{44} and~\ref{subsetsum4} in case
$g_0 \leq 4.$

\vspace{2mm}
To prove the theorem, we need to show that if $s'' > p-1$, then
$s(\sqrt{2p})\in [p+2,p+7].$
From claim (v) and Lemma~\ref{except-one}, it follows that when $s''
>p-1$, then we have $\delta(p)=1.$ From Theorem~\ref{card}, it follows that
\begin{equation}\label{ls2p}
s(\sqrt{2p}) \geq p+2.
\end{equation}
Moreover from Lemmas~\ref{subsetsum5}
and~\ref{subsetsum4}, it follows that the set $[3, s'']_p$ is contained
in $(\aaa)^{\sharp}$ in case $g_0 \notin \{1,2\}$. Therefore it
follows that if $s''> p-1$, then we have
$$g_0 \in \{1,2\}.$$
When $g_0 \in \{1,2\}$, then we have
\begin{equation}\label{gs2p}
s(\sqrt{2p}) - g_0 \leq \sum_{a \in \aaa}|a|_p= \sum_{a'\in \aaa'}|a'|_p +s''
\end{equation}
and from Lemma~\ref{subsetsum4}, it follows
that either $s'' \leq p-1$ or we have $s'' = p+g_0.$ 
We also know all the possibilities of $\aaa'$ from Lemma~\ref{44} and
claim (i). Using this
and rearranging the terms in~\eqref{gs2p}, we obtain that
when $s'' > p-1$, then we have
\begin{equation}\label{us2p}
s(\sqrt{2p}) \leq p+7.
\end{equation}
Therefore if $s''\geq p-1$ then from~\eqref{ls2p} and~\eqref{us2p}, we have $s(\sqrt{2p}) \in [p+2,p+7]$. Hence the theorem follows.
\end{proof}
The Theorem~\ref{largest} readily follows from Theorems~\ref{card} and~\ref{struc-largest}.
\section{Proof of Theorem~\ref{sl0}}
Let $\mathcal{A}$ be as in Theorem~\ref{sl0}. From the assumptions we
have
\begin{equation}
e(\aaa) := |\sqrt{2p}-\card(\aaa)| \le \psi(p)\sqrt{p} \; \text{ and } p \text{ is sufficiently large, }
\end{equation}
where  $\psi$ is a function from $[2, \infty) \text{ to } \rr^+$ which
tends to $0$ at $\infty.$ In what follows $\psi$ will denote this function.

\vspace{2mm}
From
Proposition~\ref{mainprop}, replacing $\aaa$ by $d.\aaa$ for some
non-zero element $d \in \zp$ we have
\begin{equation}\label{smallsum}
\sum_{a \in \mathcal{A}}|a|_p \leq p + O\left((e(\aaa)^{3/2}\ln
(e(\aaa) +2 )\right)
\end{equation}
and
\begin{equation}\label{few-ve}
\sum_{a \in \aaa, a < 0 }|a|_p = O\left(e(\aaa)^{3/2} \ln (e(\aaa) +2)\right).
\end{equation}
As before we find it more convenient to work with $\baaa$ than
$\aaa$.  We partition the set of natural numbers
 into the three disjoint sets $P, N$ and $G$ which are defined
as follows.
\[
P = \{ k \in \N  :               k \in \baaa\}, \ N = \{ k : k \in -\baaa\}, \ G =
\{k : k \notin \baaa \cup -\baaa\}.
\] 
An immediate corollary of~\eqref{few-ve} is that the cardinality of
$N$ is $O\left(e(\aaa)^{3/4}\ln e(\aaa)\right)$. 
We shall prove the following result.
\begin{prop}\label{nsmall}
The cardinality of $N$ is $O(\sqrt{e(\aaa)})$. Moreover there exists an
  absolute constant $c$ such that $N \subset [1, ce(\aaa)].$
\end{prop}
We first deduce Theorem~\ref{sl0} from Proposition~\ref{nsmall}.
\begin{proof}[Proof of Theorem~\ref{sl0}]
From Proposition~\ref{mainprop} there exists a $d\in (\zp)^*$ such
that replacing $\aaa$ by $d.\aaa$, the inequalities~\eqref{smallsum}
and~\eqref{few-ve} hold. Let $\kkk = \baaa \setminus (-N)$.
Then we
have $e(\kkk): = |\sqrt{2p}-\card(\aaa)| \leq e(\aaa) + \card(N) =
O\left(\psi(p) \sqrt{p}\right),$ the last equality follows using
Proposition~\ref{nsmall}.
 Moreover we have
$$\sum_{k \in \kkk} |k| \leq \sum_{\bar{a} \in \baaa}|\bar{a}| \leq p
+ O(e(\aaa)^{3/4} \ln(e(\aaa)+2)).$$
Therefore it follows that $\kkk$ satisfies the assumption of
Proposition~\ref{zprop} with $s(\kkk) =  O(e(\aaa)^{3/4}
\ln(e(\aaa)+2)).$
Let $\ccc_1$ be a subset of $\kkk$ as in the proof of
Proposition~\ref{zprop}. Then we have 
 $\mathcal{C}_1 \subset [1,
ce(\aaa)]$, $\card(\mathcal{C}_1) = O(\sqrt{e(\aaa)} \ln e(\aaa) )$
and $\sum_{k \in \kkk \setminus \mathcal{C}_1 } |k| \leq p
  -1$. Let $\aaa' = \sigma_p(N \cup \mathcal{C}_1)$ and
  $\aaa'' = \aaa \setminus \aaa'$. Then using Proposition~\ref{nsmall}
  and the properties of $\ccc_1$ just stated, we have that $\aaa' \subset [-ce(\aaa), ce(\aaa)]_p$ for some
  absolute positive constant $c$ and $\card(\aaa') = \card(\ccc_1) +
  \card(N) = O(\sqrt{e(\aaa)\ln(e(\aaa)+2)}).$ From the definition of
  $N$ and
  $\aaa''$, we have that $\aaa'' \subset [1, \frac{p}{2}]_p$. Moreover
  we have
$$\sum_{a'' \in \aaa''}|a''| = \sum_{k \in \kkk \setminus \mathcal{C}_1 } |k| \leq p
  -1.$$
 
Hence
Theorem~\ref{sl0} follows.                              
\end{proof}
\subsection{Proof of Proposition~\ref{nsmall}}
\begin{lem}\label{lotof+ve}
The cardinality of $P(0.9\sqrt{2p})$ is equal to $0.9\sqrt{2p} - O(e(\aaa))$. 
\end{lem}
\begin{proof}
Applying Lemma~\ref{dense-int} with $\kkk = \baaa$ and $e(p) =
e(\aaa)^{3/2}\ln e(\aaa)$ we obtain that 
\[
\card(P(0.9\sqrt{2p})) + \card(N(0.9\sqrt{2p})) = 0.9\sqrt{2p} - O(e(\aaa)) 
\]
and using~\eqref{few-ve} it also follows that the cardinality of $N$ is $O(e(\aaa)^{3/4}\ln
e(\aaa))$. Hence the lemma follows.
\end{proof}
\begin{lem}\label{13/8}
Let $q$ be a sufficiently large positive integer and $B \subset [1,q]$ with
$\card(B) \geq \frac{7}{8}q$. Then the interval $[q + 1,
\frac{13}{8}q]$ is contained in $2\verb|^|B$.
\end{lem}
\begin{proof}
For any $n \in [q+1, \frac{13}{8}q]$ there are $q - [\frac{n}{2}] -1$
pairs of elements $(a_i, b_i)$ with $n = a_i + b_i$, $a_i < b_i$ and
both $a_i, b_i \in [1, q]$. Among these pairs if there is a pair
$(a_i, b_i)$ with both $a_i, b_i \in B$ then the assertion follows. If
not then $\card(B) \leq q - (q -  [\frac{n}{2}] -1) = [\frac{n}{2}]
+ 1$ which is strictly less than $\frac{7}{8}q$, since $n \leq
\frac{13}{8}q$. This is contrary to the assumption. Hence the lemma follows.
\end{proof}
\begin{lem}\label{sP}
Let $q$ be a sufficiently large positive integer and $B \subset [1,q]$
with $\card(B) = q - O(\psi(q)q)$. Then the interval $[q+1,
\psi(q)^{1/2}q^2]$ is contained in the set $B^{\sharp}$.
\end{lem}
\begin{proof}
For any $n \in [q+1, \frac{13}{8}q]$ it follows from the previous
lemma that $n \in B^{\sharp}$. 
Let $B(0.2q, 0.4q) = B \cap [0.2q, 0.4q] = \{b_1 > b_2 >..... >
b_I\}$. Then from the assumptions of the lemma we have $\card(B(0.2q,
0.4q)) \geq 0.2q - O(\psi(q)q)$. Let $C$ be the sequence $\{c_i\}_{i =
1}^{I}$ with $c_i = \sum_{l =1}^{i}b_l$.  Then the following
properties of $c_i$ are evident.
\begin{enumerate}
 \item $c_{i} \geq 0.2q i$,              
\item $c_{i+1} - c_{i} \leq 0.4q$.
\end{enumerate}
Now for every $n \in [\frac{13}{8}q,
\psi(q)^{1/2}q^2]$, let $n_i$ be the least integer with $1 \leq n_i
\leq I$ such that $n-c_{n_i}$ belongs to  the interval $[1.01 q,
\frac{13}{8}q]$. From the properties of $c_i$ it follows that such a
$n_i$ exists and $n_i \leq \psi(q)^{1/2}q$. Moreover we also have $c_i
\in B_{n_i}^{\sharp}$, where $B_{n_i} = \{b_1, b_2 , \ldots, b_{n_i}\}
\subset B$ and is of cardinality $n_i$. Now $\card(B \setminus B_{n_i})
\geq q - O(\psi(q)^{1/2}q)$. Therefore using Lemma~\ref{13/8}, the
element $n- c_{n_i}$ can be written as a sum of distinct elements of
the set $B \setminus B_{n_i}$. Hence $n \in B^{\sharp}$. Hence the
lemma follows. 
\end{proof}
\begin{lem}\label{lP}
The set $P^{\sharp}$ contains the interval $[0.9\sqrt{2p}+1, \psi(p)^{1/2}p]$.
\end{lem}
\begin{proof}
From Lemma~\ref{lotof+ve}, the cardinality of $P(0.9\sqrt{2p})$ is
$0.9\sqrt{2p} - O(e(\aaa)) \geq 0.9\sqrt{2p} -
O(\psi(p)\sqrt{p})$. Therefore the assertion follows
from Lemma~\ref{sP}.
\end{proof}
\begin{lem}\label{nless}
The cardinality of $N$ is $O(\sqrt{e(\aaa)})$.
\end{lem}
\begin{proof}
From Lemma\ref{a-bdd}, the largest integer $y_0$ belonging to $P
\cup N$ is $O(e(\aaa)\sqrt{2p})$.
Since $\baaa$ does not contain any multiple of $p$ and hence does not
contain zero, the sets $P^{\sharp}$ and $N$ are disjoint. Therefore using
Lemma~\ref{lP}, it follows that $N \subset [1, 0.9\sqrt{2p}]$. Since
the cardinality of $N$ is  $O(e(\aaa)^{3/4}\ln e(\aaa)$, it follows
that $N^{\sharp} \subset [1, c_0e(\aaa)^{3/4}\ln e(\aaa) \sqrt{p}]$. Since
$e(\aaa) \leq \psi(p) \sqrt{p}$, using Lemma~\ref{lP}, it follows
that $N^{\sharp} \subset [1, 0.9\sqrt{2p}]$. Now using
Lemma~\ref{lotof+ve} and the fact that $P$ and $N^{\sharp}$ are disjoint sets,
 it follows that the cardinality of $N^{\sharp}$ is
$O(e(\aaa))$. Since we also have that the cardinality of $N^{\sharp}$ is at
least $\frac{(\card(N))^2}{2}$, the assertion follows.
\end{proof}
\begin{lem}\label{naresmall}
There exists a positive absolute constant $c_0$ such that $N \subset
[1, c_0e(\aaa)]$.
\end{lem}
\begin{proof}
Let $x$ be a sufficiently large integer such that $\card(P(x)) \geq
\frac{7}{8}x$, then using Lemma~\ref{sP}, the set $N$ does not contain
any element from the interval $[x+1, \frac{13}{8}x]$. From
Lemma~\ref{lotof+ve}, there exists an integer $x_0$ such that $x_0 =
O(e(\aaa))$ and for any integer $x$ with $x_0 \leq x  \leq
0.9\sqrt{2p}$, we have $\card(P(x)) \geq \frac{7}{8}x$. Therefore the
set $N$ does not contain any integer in the interval $[x_0,
0.9\sqrt{2p}]$. As it was observed during the proof of
Lemma~\ref{nless}, we have
$N \subset [1, 0.9\sqrt{2p}]$, it follows that $N \subset [1,
x_0]$. Hence the lemma follows.
\end{proof}
From Lemmas~\ref{nless} and \ref{naresmall}, Proposition~\ref{nsmall} follows.

\noindent
Jean-Marc Deshouillers\\
Institut Math{\'e}matique de Bordeaux,\\
Universit\'e de Bordeaux et CNRS\\
F-33405 TALENCE Cedex,\\
France.\\
E-mail: jean-marc.deshouillers@math.u-bordeaux1.fr\\

\noindent
Gyan Prakash\\
Institut Math{\'e}matique de Bordeaux,\\
Universit\'e de Bordeaux 1,\\
F-33405 TALENCE Cedex,\\
France.\\
E-mail: gyan.prakash@math.u-bordeaux1.fr\\
gyan.jp@gmail.com
\end{document}